\documentclass[12pt]{amsart}

\date{May 30, 2016}

\usepackage[margin=1.15in]{geometry}

\usepackage{amsmath,amscd,amssymb,amsfonts,latexsym,wasysym, mathrsfs, enumitem, mathtools,hhline,color}
\usepackage[all, cmtip]{xy}

\usepackage{url}

\definecolor{hot}{RGB}{65,105,225}

\usepackage[pagebackref=true,colorlinks=true, linkcolor=hot ,  citecolor=hot, urlcolor=hot]{hyperref}%make all the references and links clickable

\usepackage{graphicx,enumerate}

\newcommand{\CN}{\mathbb{C}^{n}}
\newcommand{\C}{\mathbb{C}}

\newcommand{\Z}{\mathbb{Z}}
\newcommand{\K}{\mathcal{PZ}}
\newcommand{\Su}{{\rm Supp}}
\newcommand{\V}{\mathcal{V}}
\newcommand{\m}{\mathbf{m}}
\newcommand{\G}{\mathcal{G}}

\theoremstyle{plain}
\newtheorem{theorem}{Theorem}[section]
\newtheorem{prop}[theorem]{Proposition}

\newtheorem{lm}[theorem]{Lemma}

\newtheorem{cor}[theorem]{Corollary}
\newtheorem{conj}[theorem]{Conjecture}

\newtheorem{thrm}[theorem]{Theorem}

\theoremstyle{definition}

\newtheorem{defn}[theorem]{Definition}

\newtheorem{rmk}[theorem]{Remark}

\newtheorem{ex}[theorem]{Example}
\newtheorem*{ex*}{Example}

\def\be{\begin{equation}}
\def\ee{\end{equation}}

\def\bt{\begin{thrm}}
\def\et{\end{thrm}}

\def\bc{\begin{cor}}
\def\ec{\end{cor}}

\def\br{\begin{rmk}}
\def\er{\end{rmk}}

\def\bp{\begin{prop}}
\def\ep{\end{prop}}

\def\bl{\begin{lm}}
\def\el{\end{lm}}

\def\bex{\begin{ex}}
\def\eex{\end{ex}}

\def\bd{\begin{defn}}
\def\ed{\end{defn}}

\newcommand\sS{\mathcal{S}}

\def\cZ{\mathcal{Z}}

\DeclareMathOperator{\Supp}{Supp}                % Supp
\DeclareMathOperator{\Exp}{Exp}

\def\ra{\rightarrow}

\def\bC{\mathbb{C}}

\def\cV{\mathcal{V}}

\def\al{\alpha}

\def\bP{\mathbb{P}}
\def\cH{\mathcal{H}}

\def\cD{\mathcal{D}}
\def\cO{\mathcal{O}}
\def\lra{\longrightarrow}

\def\lam{\lambda}

\def\bN{\mathbb{N}}

\newcommand{\ubul}{{\,\begin{picture}(-1,1)(-1,-3)\circle*{2}\end{picture}\ }}

\title{Cohomology support loci of local systems}

\author{Nero Budur}
\author{Yongqiang Liu}
\author{Luis Saumell}
\author{Botong Wang}

\address{KU Leuven, Department of Mathematics,
Celestijnenlaan 200B, B-3001 Leuven, Belgium} 

\email{nero.budur@kuleuven.be}
\email{yongqiang.liu@kuleuven.be}
\email{luis.saumell@kuleuven.be} 
\email{botong.wang@kuleuven.be}



\begin{document}

\maketitle

\begin{abstract} The support S of Sabbah's specialization complex is a simultaneous generalization of the set of eigenvalues of the monodromy on Deligne's nearby cycles complex, of the support of the Alexander modules of an algebraic knot, and of certain cohomology support loci. Moreover, it equals conjecturally the image under the exponential map of the zero locus of the Bernstein-Sato ideal. Sabbah showed that S is contained in a union of translated subtori of codimension one in a complex affine torus. Budur-Wang showed recently that S is a union of torsion-translated subtori. We show here that S is always a hypersurface, and that it admits a formula in terms of log resolutions. As an application, we give a criterion in terms of log resolutions for the (semi-)simplicity as perverse sheaves, or as regular holonomic D-modules, of the direct images of rank one local systems under an open embedding. For hyperplane arrangements, this criterion is combinatorial.
 \end{abstract}

\section{Introduction}

\subsection{Cohomology support loci.} For a topological space $T$, let $M_B(T)$ be the moduli space of rank $1$ $\bC$-local systems on $T$. The {\it cohomology support loci of $T$} are defined as $$\V(T)=\lbrace L\in M_{B}(T) \mid \dim H^\ubul(T,L) \ne 0 \rbrace,$$
and are homotopy invariants of $T$. It was shown recently in \cite{BW-a,BW} that $\V(T)$ are finite unions of torsion translated affine subtori of the affine algebraic group $M_B(T)\cong Hom(H_{1}(T, \Z),\C^{\ast})$, if $T$ is a smooth complex quasi-projective algebraic variety, or a small ball complement of the germ of a complex analytic set in a complex manifold. It remains though a difficult task to compute cohomology support loci. This article is an application of the structure result for cohomology support loci.

Let $j:U\ra X$ be the open embedding in a complex manifold $X$ of the complement of a hypersurface $f^{-1}(0)$, where $f:X\ra \bC$ is an non-invertible analytic function. For $x\in f^{-1}(0)$, let $U_x$ be the complement in a small ball in $X$ centered at $x$ of  $f^{-1}(0)$. Then $M_{B}(U_{x})\cong (\bC^*)^r$,
where $r$ is the number of analytic branches of $f$ at $x$. In this article, we show that if we take the union, in a certain sense, of the cohomology support loci $\V(U_x)$ for all points $x\in f^{-1}(0)$, the resulting set $\V(U,X)$ is much easier to deal with and that it contains a wealth of geometric information.

More precisely, for a point $x\in f^{-1}(0)$, define the map
$$
res_x:M_B(U)\ra M_B(U_x)
$$ 
as the restriction of a local system on $U$ to $U_x$. Define 
$$
\V(U,X) =\mathop{\bigcup_{x\in f^{-1}(0)}} res_x^{-1}(\V(U_x)) \quad\subset M_B(U).
$$
Note that with these general assumptions on $X$ and $U$, it can happen that $M_B(U)$ is infinite-dimensional and $\V(U,X)$ has infinitely many irreducible components. Note that $\V(U,X)$ can be similarly defined for any open proper subset $U$ of a complex analytic space $X$, namely, as the the set of rank-1 local systems on $U$ with non-trivial cohomology locally on $U$ along the complement $X\setminus U$.

We show that $\V(U,X)$ admits a simple formula in terms of log resolutions, in a sense to be made clear below. Consider first the case $U=U_x$ and $X=B_x$ is a small ball at $x$. That is $f$ is the germ of an analytic function. We can assume that $f$ is reduced and that $F=(f_1,\ldots, f_r)$ is the collection of reduced analytic branches of $f$. In particular $f=\prod_{i=1}f_i$ on $B_x$. Let $\mu:Y\ra X$ be a log resolution of $f$ that also blows up the point $x$. Let $E_j$ with $j\in J$ be the irreducible components of $(f\circ\mu)^{-1}(0)$. Let $a_{ij}$ be the order of vanishing of $f_i$ along $E_j$. Let $E_j^\circ=E_j\setminus \cup_{i\ne j}(E_i\cap E_j)$. 
Define 
$$
Z^{mon}_{F,x}(t_1,\ldots,t_r)= \mathop{\prod_{j\in J\text{ with }}}_{\mu(E_j)=x}\left ( t_1^{a_{1j}}\cdot\ldots\cdot t_r^{a_{rj}} -1 \right )^{-\chi(E_j^\circ)}.
$$
By \cite{Sab}, this is the {\it multi-variable monodromy zeta function at $x$} of Sabbah's specialization complex, which we will bring into focus soon.  It is also a multi-variable version of a classical formula of A'Campo \cite{AC}.

For a rational function $Q(t_1,\ldots, t_r)$, let
$
\cZ(Q)\text{ and }\K(Q)
$
denote the zero locus and, respectively, the union of the polar and the zero locus of $Q$.

With this notation, we prove:

\bt\label{thrm1}
Let $F=(f_1,\ldots, f_r):(X,x)\ra (\bC^r,0)$ be a collection of germs of non-invertible irreducible analytic functions on a complex manifold. Let $f=\prod_{i=1}^rf_i$. Assume that for all points $y\in f^{-1}(0)$ close to $x$, the germ of $f$ at $y$ is reduced. Then
$$
\V(U_x,B_x) = \mathop{\bigcup_{y\in f^{-1}(0)}}_{y\text{ close to }x} \K (Z^{mon}_{F,y})\quad\subset M_B(U_x)=(\bC^*)^r,
$$
where the union is over generic points $y$ of the finitely many strata of a Whitney stratification for $f^{-1}(0)$. In particular, $\V(U_x,B_x)$ is a finite union of torsion translated codimension-one subtori of $(\bC^*)^r$.
\et

It is well-known that $\V(U_x,B_x)$ is a generalization of the set of eigenvalues of the monodromy of the cohomology of the Milnor fibers of $f$ at points ranging over its zero locus, see \cite[Proposition 1.3]{BW}. Therefore this theorem generalizes an observation of J. Denef \cite{D} that every such eigenvalue for $f$ appears as a pole or a zero of the monodromy zeta function of $f$ at some point $y\in f^{-1}(0)$.

As a corollary, one obtains a simple combinatorial formula in the case of hyperplane arrangements:

\bt\label{thrmHA}
Let $F=(f_1,\ldots,f_r)$ be a collection of linear forms on $X=\bC^n$ defining mutually distinct hyperplanes. Then
$$
\V(U,X) = \cZ \left(\prod_W\left ( \prod_{i\; :\; f_i(W)=0}t_i-1\right)\right),
$$
where the first product is over the dense edges $W$ of $f=\prod_if_i$.
\et

Theorem \ref{thrm1} is a special case of more general result, as we explain next.

\subsection{Sabbah's specialization complex.} Sabbah's specialization complex, introduced in \cite{Sab}, is a simultaneous generalization of  Deligne's nearby cycles complex and of the Alexander modules of an algebraic knot. 
Given a collection $F=(f_1,\ldots,f_r)$ of analytic functions $f_i:X\ra \bC$ on a complex manifold, Sabbah  defined a complex $\psi_F(\bC_X)$ with $A$-constructible cohomology on $f^{-1}(0)$, where $A=\bC[t_1^{\pm},\ldots,t_r^{\pm}]$ is the affine coordinate ring of $(\bC^*)^r$, and $f=\prod_if_i$. This complex is the analog of Deligne's nearby cycles complex $\psi_f(\bC_X)$ for the case $r=1$. While $\psi_f(\bC_X)$ governs the Milnor monodromy information, Sabbah's $\psi_F(\bC_X)$ governs  the more general Alexander-type invariants. 

One of the main results of \cite{Sab} is about the support in $(\bC^*)^r$ of the stalks of $\psi_F(\bC_X)$ given by the $A$-module structure, denoted 
$$
\Supp_x(\psi_F(\bC_X)).
$$
$
\Supp_x(\psi_F(\bC_X))
$ is shown in {\it loc. cit.} to be included in a special hypersurface whose irreducible components are translated subtori.

In \cite{BW}, it was shown that each component of $\Supp_x(\psi_F(\bC_X))$ is a torsion translated subtorus of $(\bC^*)^r$. This follows from the relation with cohomology support loci, as we describe now. Identify $M_B((\bC^*)^r)$ with $(\bC^*)^r$ via monodromies around the coordinate axes. Let
$$
\gamma_{U,F}: (\bC^*)^r=M_B((\bC^*)^r)\ra M_B(U)
$$
be the map which pulls back local systems from $(\bC^*)^r$ to $U$ via $F$, where $U=X\setminus f^{-1}(0)$.
Then, by \cite{B, LM},  \be\label{eS}
\Supp_x(\psi_F(\bC_X))
=\gamma_{U,F}^{-1} (res_x^{-1}(\V(U_x))). \ee
Define the {\it support of Sabbah's specialization complex} to be 
$$\sS(F)=\mathop{\bigcup_{x\in f^{-1}(0)}}\Supp_x(\psi_F(\bC_X))\mathop{=\joinrel=}^{(1)} \gamma_{U,F}^{-1}(\V(U,X))\quad\subset (\bC^*)^r.$$
Equivalently, easier to remember but less precise, $$\sS(F)= \{\text{rank one local systems on }(\bC^*)^r\text{ with }H^\ubul \ne 0\text{ locally on }U$$
$$\text{ along }f^{-1}(0)\text{ under the map }F:U\ra (\bC^*)^r\}.$$

Theorem \ref{thrm1} is then a particular case of the following surprisingly strong improvement of the original result of Sabbah:

\bt\label{thrmSF}
Let $F=(f_1,\ldots, f_r):X\ra \bC^r$ be a collection of non-invertible analytic functions on a complex manifold. Let $f=\prod_{i=1}^rf_i$ and assume that $f^{-1}(0)$ admits a finite Whitney stratification. Then
$$
\sS(F) = \mathop{\bigcup_{x\in f^{-1}(0)}} \K (Z^{mon}_{F,x}),
$$
where the union is over generic points $x$ of the Whitney strata. In particular, $\sS(F)$ is a finite union of torsion translated codimension-one subtori of $(\bC^*)^r$.
\et

We show in Example \ref{exC1} that  $\V(U,X)$ can fail to be of pure codimension 1. Without the finiteness assumption on the number of Whitney strata, $\sS(F)$ could be a countable union of torsion-translated codimension-one subtori. This assumption is satisfied for example in the algebraic case, or in the local analytic case.

\subsection{Bernstein-Sato ideals}
Let us mention that a different way to compute $\sS(F)$ without appealing to a log resolution was conjectured in \cite{B}: 

\begin{conj} Let $F=(f_1,\ldots, f_r):X\ra \bC^r$ be a collection of non-invertible analytic functions on a complex manifold such that $f^{-1}(0)$ admits a finite Whitney stratification, where $f=\prod_{i=1}^rf_i$. Then $$\sS(F)=\Exp (\cZ(B_F)),$$ where $\Exp:\bC^r\ra(\bC^*)^r$ is the map $\al\mapsto \exp(2\pi i\al)$, and $\cZ(B_F)$ is the zero locus of the Bernstein-Sato ideal of $F$. \end{conj}

It is also conjectured in {\it loc. cit.} that the Bernstein-Sato ideal $B_F$, although not necessarily principal, is generated by products of linear polynomials of type $\sum_{i=1}^rb_is_i + b$ with $b_i\in\bN$, $b\in\bN\setminus\{0\}$. Theorem \ref{thrmSF} suggests surprisingly that more might be true, namely that the zero locus $\cZ(B_F)$ is of pure codimension-one. All these conjectures are of course true for the $r=1$ case, by classical results of Malgrange and Kashiwara.

It was proved in \cite{B} that one inclusion holds for the above conjecture:
\be\label{cp}
\sS(F)\subset\Exp (\cZ(B_F)).
\ee
We fix in Remark \ref{remBF} a gap in the proof in {\it loc. cit.} of this inclusion.

\subsection{Semi-simplicity of direct images.} Let $j:U\ra X$ be the open embedding of the complement of a non-empty hypersurface in a complex manifold of dimension $n$. For a local system $L$ of rank one on $U$, the shifted complex $L[n]$ is perverse on $U$. The derived direct image $Rj_*(L[n])$ and the direct image with compact supports $j_!(L[n])$ are perverse sheaves on $X$. Let $Perv(X)$ denote the abelian and artinian category of perverse sheaves on $X$. We address the question of how to detect if $Rj_*(L[n])$ and $j_!(L[n])$ are semi-simple. In this case, simplicity and semi-simplicity are equivalent.

In the space of rank one local systems on $U$, $M_B(U)$, define the {\it non-simple locus}
$$
\V^{ns}(U,X)=\{ L\in M_B(U)\mid Rj_*(L[n])\neq \text{ (semi-)simple in }Perv(X) \}.
$$
We prove the following (semi)-simplicity criterion and relation with cohomology support loci:

\bt\label{thrmNS} Let $U$ be the complement of a non-empty hypersurface in a complex manifold $X$. Then 
$$\V^{ns}(U,X)=\cV(U,X).$$
\et

Moreover, since $\V(U,X)$ is stable under taking the inverses of local systems (see for example \cite[Theorem 1.2]{BW} for a more general statement), $\V(U,X)$ is also the locus of local systems $L$ with $j_!(L[n])$ not (semi-)simple.  

In particular:

\bc\label{corFns} If $F=(f_1,\ldots,f_r):X\ra\bC^r$ is a collection of non-invertible analytic functions, $f=\prod_{i=1}^rf_i$, $U=X\setminus f^{-1}(0)$, then 
$$\sS(F)=\gamma_{U,F}^{-1}(\V^{ns}(U,X)).$$
\ec

Hence the (semi-)simplicity question is answered in terms of log resolutions, and conjecturally from the Bernstein-Sato ideal of $F$ if $f^{-1}(0)$ admits only finitely many Whitney strata.

Since perversity is a local condition, the theorem follows from the local case, namely from 
$$\V^{ns}(U_x,B_x)=\cV(U_x,B_x),$$
for all $x\in X\setminus U$.

By the Riemann-Hilbert correspondence, the (semi-)simplicity criterion has a $\cD$-module counterpart. Let $\cD_{X}$ denote the sheaf of analytic linear differential operators on $X$. Let $Mod_{rh}(\cD_{X})$ be the category of regular holonomic left $\cD_{X}$-modules. Let $DR_X:Mod_{rh}(\cD_{X})\ra Perv(X)$ be the de Rham functor, an equivalence of categories. Let $F$ be as in Corollary \ref{corFns}. For $\lam\in(\bC^*)^r$, let $\al\in \Exp^{-1}(\lam)\subset \bC^r$ and $k\in\bN$, $k\gg 0$. One has the known isomorphisms in $Mod_{rh}(\cD_{X})$:
\begin{align*}
 j_*(\cD_U\prod_{i=1}^rf_i^{\al_i})&\cong\cO_{X}[f^{-1}]\prod_{i=1}^rf_i^{\al_i}\cong \cD_{X}\prod_{i=1}^rf_i^{\al_i-k}\\
&\cong \cD_{X}[s]f^s\prod_{i=1}^rf_i^{\al_i}/ (s+k)\cD_{X}[s]f^s\prod_{i=1}^rf_i^{\al_i},\\
j_{!}(\cD_U\prod_{i=1}^rf_i^{\al_i})&\cong \cD_{X}[s]f^s\prod_{i=1}^rf_i^{\al_i}/ (s-k)\cD_{X}[s]f^s\prod_{i=1}^rf_i^{\al_i}.
\end{align*}
Here $j_*$ and $j_{!}$ are the direct image and, respectively, the special direct image for regular holonomic $\cD$-modules, such that 
\begin{align*}
DR_U(\cD_U\prod_{i=1}^rf_i^{\al_i}))&\cong L_\lam[n],\\
DR_X(j_*(\cD_U\prod_{i=1}^rf_i^{\al_i}))&=Rj_*(L_\lam[n]),\\
DR_X(j_{!}(\cD_U\prod_{i=1}^rf_i^{\al_i}))&=j_!(L_\lam[n]).
\end{align*}

Now one can easily make the translation of the (semi-)simplicity criterion from the above theorem into $\cD$-modules:

\bc\label{corBF} With the assumptions as in Corollary \ref{corFns},
\begin{align*}
\sS(F)&=\{ \lam\in (\bC^*)^r\mid j_*(\cD_U\prod_{i=1}^rf_i^{\al_i})\neq\text{ (semi-)simple in } Mod_{rh}(\cD_X) \}\\
&=\{ \lam\in (\bC^*)^r\mid j_!(\cD_U\prod_{i=1}^rf_i^{\al_i})\neq\text{ (semi-)simple in } Mod_{rh}(\cD_X) \}.
\end{align*}
\ec
 
Together with Theorem \ref{thrmHA} and Theorem \ref{thrmNS}, this gives a very easy necessary and sufficient combinatorial criterion for simplicity of direct images in the case of hyperplane arrangements. A limited sufficient criterion had been obtained earlier by Abebaw-B\o{}gvad \cite{AB}.

\subsection{Acknowledgement.} The authors were partly sponsored by a FWO grant, a KU Leuven OT grant, and a Flemish Methusalem grant.

\section{Sabbah's specialization complex}

In this section we recall a few facts about Sabbah's specialization complex from \cite{Sab, B, LM}.

Let $F=(f_1,\ldots,f_r):X\ra \bC^r$ be a collection of analytic functions on a complex manifold $X$ of dimension $n$. Let $f=\prod_{i=1}^rf_i$ and $U=X\setminus f^{-1}(0)$. Let $j:U\ra X$ and $i:f^{-1}(0)\ra X$ be the natural open and, respectively, closed embeddings.

Consider the following commutative diagram of fibered squares of natural maps:\begin{center}
$\xymatrix{
f^{-1}(0) \ar@{^{(}->}[r]^{i}   & X  \ar[d]^{F}   & \ar@{_{(}->}[l]_{j} U \ar[d]^{F}  & \ar[l]_{\pi} \widehat{U} \ar[d]^{\widehat{F}} \\
                       & \C^{r}    & \ar@{_{(}->}[l] (\C^{\ast})^{r}                                           
            & \ar[l]_{\widehat{\pi}} \widehat{(\C^{\ast})^{r}} 
}$
\end{center}
where $\widehat{\pi}=\Exp:\widehat{(\C^{\ast})^{r}}\ra (\C^*)^r$ is the universal covering.

\bd  \label{d2.1}  {\it Sabbah's specialization complex functor}  of $F$ is defined by \begin{center}
$\psi_{F}=i^{-1}Rj_{\ast} R\pi_{!} (j \circ \pi)^{\ast} : D^{b}_{c}(X,\C) \rightarrow D^{b}_{c}(f^{-1}(0),A)$, 
\end{center} where $A=\C[t_{1},t_{1}^{-1},\cdots,t_{r},t_{r}^{-1}]$. Here $D^{b}_{c}(X,\C)$ is the derived category of bounded complexes of sheaves with $\bC$-constructible cohomology on $X$, and $D^{b}_{c}(f^{-1}(0),A)$ is the derived category of bounded complexes of sheaves with $A$-constructible cohomology on $f^{-1}(0)$. We call $\psi_{F}(\C_{X})$  {\it Sabbah's specialization complex.} 
\ed

\bl\label{lemBr}{\rm (\cite{Br})} When $r=1$, $\psi_f(\bC_X)$ as defined here equals the shift by $[-1]$ of Deligne's nearby cycles complex together with the action of the monodromy.
\el

 \bd \label{d2.2} For any finitely generated $A$-module $M$, the {\it support  of $M$} is the zero locus of the annihilator ideal of $M$: $$\Su(M)=\cZ(ann(M)) \subset (\C^{\ast})^{r}.$$
 \ed
 
Let $P$ be a prime ideal in $A$, of height $1$. Then $P$ is principal, and we let $\Delta({P})$ denote the generator of $P$, which is well-defined up to multiplication by units of $A$. Denote by $A_{P}$ the localization of $A$ at $P$, then $A_{P}$ is a principal ideal domain.   Assume that $\Su (M)$ is proper in $(\C^{\ast})^{r}$, that is, $\Su (M)$ has codimension  at least 1 in $(\C^{\ast})^{r}$. Then $M_{P}$ has finite length as a $A_{P}$-module, which is denoted by $lg(M_{P})$. The {\it characteristic polynomial of $M$} is defined as $$ \Delta(M)=\prod_{P} \Delta({P})^{lg(M_{P})},$$ 
where the product is over all the prime ideals  in $A$ of height 1 such that $\cZ(P)\subset \Su(M)$. Since $\Su (M)$ is proper, this product is indeed a finite product. The prime factors of $\Delta(M)$ are in one-to-one correspondence with the codimension one irreducible hypersurfaces of $(\C^{\ast})^{r}$ contained in $\Su(M)$.

\bd \label{d2.3}  For $\G \in D_{c}^{b}(X,A)$ and a point $x\in X$, the  {\it support of $\G$ at $x$} is defined by
$$\Su_{x} (\G) := \bigcup_{i} \Su (\mathcal{H} ^{i}(\G)_{x}) \subset (\C^{\ast})^{r},$$ and the {\it multi-variable monodromy zeta-function of $\G$ at $x$} is defined by
$$Z_{x}^{mon}(\G)(t_1,\ldots,t_r):=\prod_{i} \Delta (\mathcal{H} ^{i}(\G)_{x})^{(-1)^{i}}\quad\in \bC(t_1,\ldots,t_r).$$
\ed
 
 \bd \label{d2.4} 
 The {\it support of Sabbah's specialization complex} is defined to be 
$$\sS(F)=\mathop{\bigcup_{x\in f^{-1}(0)}}\Supp_x(\psi_F(\bC_X)).$$
 The {\it multi-variable monodromy zeta function of  $F$ at $x\in f^{-1}(0)$} is defined as 
 $$  Z_{F,x}^{mon}(t_1,\ldots,t_r)= Z_{x}^{mon}(\psi_{F} (\C_X))(t_1,\ldots,t_r).$$
 \ed

Let $\mu:Y\ra X$ be a log resolution of $f$ that also blows up the point $x$. Let $E_j$ with $j\in J$ be the irreducible components of $(f\circ\mu)^{-1}(0)$. Let $a_{ij}$ be the order of vanishing of $f_i$ along $E_j$. Let $E_j^\circ=E_j\setminus \cup_{i\ne j}(E_i\cap E_j)$. 
Then, we have the following generalization of A'Campo's formula: 
\bt {\rm (\cite[2.6]{Sab})}
$$
Z^{mon}_{F,x}(t_1,\ldots,t_r)= \mathop{\prod_{j\in J\text{ with }}}_{\mu(E_j)=x}\left ( t_1^{a_{1j}}\cdot\ldots\cdot t_r^{a_{rj}} -1 \right )^{-\chi(E_j^\circ)}.
$$
\et

In the introduction, we have used the following:

\bl {\rm{(\cite{B, LM})}} Keeping the notation from above and from the introduction:

(a) $\Supp_x(\psi_F(\bC_X))
=\gamma_{U,F}^{-1} (res_x^{-1}(\V(U_x))),$ and

(b) $\sS(F)=\gamma_{U,F}^{-1}(\V(U,X))$.
 
\el 
 
\br Here, $res_y^{-1}(\V(U_y))$ is what was called in \cite{B, LM} the uniformization $\V(U_y)^{unif}$ of $\V(U_y)$  with respect to $M_B(U_x)$, for $y\in f^{-1}(0)$ close to $x$. As pointed out in \cite{LM}, in all the statements in \cite{B} where the uniform support $Supp^{unif}_x(\psi_F(\bC_X))$  appears, the {\it unif} should be dropped to conform to what is proven in \cite{B}. Indeed, the support $\Supp _x(\psi_F (\bC_X))$ needs no uniformization.
\er

 \begin{thrm}\label{thrmSab} {\rm{(\cite[Theorem 1.4]{BW})}} With the notation as above, every irreducible component of $\Supp_x(\psi_F(\bC_X))$ is a torsion translated subtorus of $(\bC^*)^r$ for all $x\in f^{-1}(0)$.
\end{thrm}

\bd Let $M=(m_{kj})\in \bN^{p\times r}$. The {\it specialization} $F^M$ of $F$ by $M$ is the map $F^M:X\ra \bC^p$ given by  $$x\mapsto (f_1^{m_{11}}\ldots f_r^{m_{1r}}(x),\ldots,f_1^{m_{p1}}\ldots f_r^{m_{pr}}(x)).$$
The specialization $F^M$ is {\it non-degenerate} if the induced map on tori $(\bC^*)^r\ra (\bC^*)^p$ given by $M$ is surjective and $\sum_{k=1}^pm_{ki}\ne 0$ for all $i$ such that $f_i$ is non-invertible.
\ed

For the next theorem, see \cite[Proposition 3.31]{B}, or \cite[2.3.8]{Sab}:

\bt \label{thrmSp}
If $G=F^M$ is a non-degenerate specialization of $F$, then for all $x$
$$
\tau_M^{-1}(\Supp_x(\psi_F(\bC_X)))=\Supp_x(\psi_G(\bC_X)),
$$
where 
\begin{gather*}
\tau_M :(\bC^*)^p\ra (\bC^*)^r\\
\tau_M:(\lam_1,\ldots,\lam_p)\mapsto (\lam_1^{m_{p1}}\ldots \lam_p^{m_{p1}},\ldots,\lam_1^{m_{1r}}\ldots\lam_p^{m_{pr}}).
\end{gather*}
\et

\section{$\sS(F)$}

\subsection{Proof of Theorem \ref{thrmSF}}

When $r=1$, the result is due to J. Denef:

\bt \label{l1}(\cite[Lemma 4.6]{D}) Let $f: (\CN,x) \to (\C,0)$ be a germ of an analytic function. If $\lambda$ is an eigenvalue of the monodromy on the cohomology of the Milnor fiber of $f$ at $x$, then  there exists  a point $y\in f^{-1}(0)$ near $x$, such that  $\lambda$ is a zero or pole of the monodromy zeta function of $f$ at $y$.
\et

We prove now the case $r>1$. We will use the notation as in the statement of Theorem \ref{thrmSF}. The inclusion of the right-hand side into the left-hand side follows from the fact that $\K (Z_{F,x}^{mon})\subset \Supp_{x}(\psi_F(\bC_X))$, by the definition of the monodromy zeta function. We prove now the reverse inclusion. 

By Theorem \ref{thrmSab}, $\Supp_{x}(\psi_F(\bC_X))$ and $\K(Z^{mon}_{F,x})$  have dense subsets of torsion points. Hence we only need to show that each torsion point of  $\sS(F)$ is contained in $\K(Z^{mon}_{F,x})$ for some $x$.

It is clear that the trivial point $(1,\cdots,1)$ is contained in $\cup_{x\in f^{-1}(0)}\K(Z^{mon}_{F,x})$.

Choose any non-trivial torsion point $P$ of $\sS(F)$. Write, by reindexing if necessary,  
$$P=\left(e^{ \tfrac{2\pi \sqrt{-1}m_{1}}{m}}, \cdots,e^{ \tfrac{2\pi \sqrt{-1}m_{l}}{m}}, 1 ,\cdots,1 \right),$$
where $m_{i}$ and $ m$ are all positive integers, and $m_{i}<m$. Without loss of generality, assume that $P\in \Supp_x(\psi_F(\bC_X))$. Set  $m_{i}=m$ for $l<i\leq r$ and $\m=(m_{1},\cdots,m_{r})$. Then we have a map associated to the $r$-tuple of numbers $\m$: $$\tau_{\m} : \C^{\ast} \to (\C^{\ast})^{r},$$ given by $ \tau \mapsto (\tau^{m_{1}},\cdots,\tau^{m_{r}})$.

Consider the specialization $f^{\m}=\prod^{r}_{i=1} f_{i}^{m_{i}}$ of $F$. This is a non-degenerate specialization, and hence Theorem \ref{thrmSp} applies. In particular,
$$\tau_{\m}^{-1}(\Supp_x(\psi_F(\bC_X)))=\Supp_x(\psi_{f^{\m}}(\bC_X)).$$
By Lemma \ref{lemBr}, the right-hand side of this equation is the set of the eigenvalues of the monodromy on the cohomology of the Milnor fibre of $f^{\m}$ at $x$. 

Set $\lambda=\exp(2\pi\sqrt{-1}/m)$, so that $P=(\lambda^{m_{1}},\ldots, \lam^{m_r})$. Hence $\lambda \in \tau_{\m}^{-1}(\Supp_x(\psi_F(\bC_X)))$. Theorem \ref{l1} gives that there exists  a point $y\in f^{-1}(0)$ near $x$, such that  $\lambda$ is a zero or pole of the monodromy zeta function of $f^{\m}$ at $y$. Without loss of generality, we can assume that $\lambda$ is a zero of $Z_{f^\m,y}^{mon}$.
A result of Sabbah \cite[Theorem 2.5.7 (c)]{Sab} states that the monodromy zeta function of the specialization is the specialization of the monodromy zeta function. More precisely, 
\be  Z_{f^\m,y}^{mon}(t) = Z_{F,y}^{mon}(t_{1}^{m_{1}},\cdots,t_{r}^{m_{r}}).
\ee
It follows that $P=(\lambda^{m_{1}},\cdots,\lambda^{m_{r}})$ is also a zero of $Z_{F,y}^{mon}$.  This finishes the proof of Theorem \ref{thrmSF}.

\begin{ex}\label{exC1}
While $\sS(F)=\gamma_{U,F}^{-1}(\V(U,X))$ is of pure codimension one, it is not always the case that the same is true about $\V(U,X)$. Take $X=\bP^1\times\bP^1$ and $U=\bC^*\times\bC$. Then $$X\setminus U=(\{0\}\times\bP^1)\cup(\{\infty\}\times\bP^1)\cup(\bP^1\times\{\infty\})$$
is a hypersurface in $X$. Then $M_B(U)=\bC^*$ by identifying a local system of rank one on $U$ with the monodromy around $\{0\}\times\bP^1$. Let $x=\{0\}\times\{\infty\}\in X\setminus U$. Then the pair $(U_x,B_x)$ has the homotopy type of the pair $((\bC^*)^2,\bC^2)$. Hence $M_B(U_x)=(\bC^*)^2$ with local systems identified with the monodromies around the germs of $\{0\}\times\bP^1$ and, respectively, $\bP^1\times\{\infty\}$ at $x$. In particular, the restriction map
$$
res_x:M_B(U)=\bC^*\ra M_B(U_x)=(\bC^*)^2
$$ 
is just the inclusion of the first coordinate subtorus $\bC^*\times\{1\}\hookrightarrow \bC^*\times\bC^*$. Moreover, the image of $res_x$ is contained in
$\V(U_x,B_x)=(\bC^*\times\{1\})\cup (\{1\}\times\bC^*)$. In particular, $res^{-1}_x(\V(U_x,B_x))=M_B(U)$. So also $\V(U,X)=M_B(U)$.
\end{ex}

\subsection{Proof of Theorem \ref{thrmHA}.} We refer to \cite{B} for the definition of the terminology ``dense edges". It was shown in \cite[Proposition 6.7]{B} that the right-hand side of the equation is the codimension one part of $\V(U,X)$, hence the theorem follows.

\section{Proof of Theorem \ref{thrmNS}} 
As mentioned in the introduction, it is enough to prove the local case. Namely, we will prove that
$$\V^{ns}(U_x,B_x)=\cV(U_x,B_x),$$
for all $x\in X\setminus U$. We can assume that $X\setminus U=f^{-1}(0)$ for a reduced germ of an analytic function $f:(B_x,x)\ra (\bC,0)$. 

In order to simplify the notation we will use $(X,U)=(B_x,U_x)$ from now where convenient. Let $j:U\ra X$ and $i:f^{-1}(0)\ra X$ be the natural open and, respectively, closed embedding. 

The following facts are well-known, see \cite{Di, DM}. In $D^b_c(X,\bC)$, one has the distinguished triangle
$$
j_!j^{-1}\ra id \ra i_*i^{-1}\mathop{\lra}^{+1},
$$
where $id$ is the identity functor. Applying this to $Rj_*(L[n])$, where $L$ is a local system of rank one on $U$, one obtains a distinguished triangle
$$
j_!(L[n])\ra Rj_*(L[n])\ra i_*i^{-1}Rj_*(L[n])\mathop{\lra}^{+1}.
$$
Since $j$ is a Stein morphism, the first two complexes on the left are actually perverse on $X$. More precisely, $j_!(L[n])$ is up to a shift the extension by zero over a closed analytic complement of the local system $L$. Hence it has constructible cohomology sheaves. Since $j$ is a Stein morphism, by \cite[5.2.17]{Di} and the remark thereafter, $j_!(L[n])$ is a perverse sheaf. Since Verdier duality preserves perversity, $Rj_*(L[n])$ is also perverse, being the Verdier dual of $j_!(L^\vee [n])$.

Hence, applying the long exact sequence of perverse cohomology, we obtain an exact sequence in $Perv(X)$,
\be\label{eqP}
0\ra{}^p\cH^{-1}(i_*i^{-1}Rj_*(L[n]))\ra j_!(L[n])\ra Rj_*(L[n])\ra {}^p\cH^{0}(i_*i^{-1}Rj_*(L[n]))\ra 0,
\ee
and the vanishing
$$
{}^p\cH^{j}(i_*i^{-1}Rj_*(L[n]))=0\quad\quad\text{ for }j\ne -1,0.
$$
By definition, the image in $Perv(X)$ of the middle map is the intersection complex of $L$, also known as the intermediate extension of the perverse sheaf $L[n]$:
$$
IC_X(L)=j_{!*}(L[n])=im\{ j_!(L[n])\ra Rj_*(L[n])\}.
$$
It is known that $IC_X(L)$ is a simple perverse sheaf on $X$. Moreover, all simple perverse sheaves on $X$ are of the type $(i_{\bar{Z}} )_* (j_Z)_{!*}(M[d_Z])$, where $Z$ is an irreducible locally closed smooth subvariety of $X$, $j_Z$ is the open embedding of $Z\ra \bar{Z}$ in its closure, $i_{\bar{Z}}$ is the inclusion $\bar{Z}\ra X$, $d_Z$ is the complex dimension of $Z$, and $M$ is an irreducible local system on $Z$.

\bl With the notation as above, let $L$ be a rank one local system on $U$. Then $Rj_*(L[n])$ is semi-simple in $Perv(X)$ if and only if it is simple. The same holds for $j_!(L[n])$.
\el
\begin{proof}
If $Rj_*(L[n])$ is semi-simple in $Perv(X)$, then one of the simple factors must be $IC_X(L)$, and the other factors must be supported on $f^{-1}(0)$. 
Taking the Verdier dual of $Rj_*(L[n])$, which is $j_!(L^{-1}[n])$, where $L^{-1}$ is the dual local system, one has then that $j_!(L^{-1}[n])$ is a direct sum of $IC_X(L^{-1})$ with factors supported on $f^{-1}(0)$. Since $j_!(L[n])=(j_!L)[n]$ and $j_!L$ is the extension by zero in this case, it follows that there are no factors supported on $f^{-1}(0)$. Hence $j_!(L^{-1}[n])$, and thus $Rj_*(L[n])$, are simple.
\end{proof}

From now on, we focus on the question about the simplicity of $Rj_*(L[n])$. If $Rj_*(L[n])$ is not simple, then in $Perv(X)$,
$$
IC_X(L)\subsetneq Rj_*(L[n]),
$$
where by the strict inclusion we mean a monomorphism that is not an isomorphism. In particular, $Rj_*(L[n])$ is not quasi-isomorphic to $j_{!}(L[n])$ in $D^b_c(X,\bC)$. The same is then true for the stalk at some point $y$ close to $x$. Hence $H^\ubul(Rj_*(L[n])_y)$ is not all zero. But
$$
H^\ubul(Rj_*(L[n])_y)=\Gamma(B_y,R^{\ubul+n} j_*(L))=H^{\ubul+n}(U_y,L).
$$
Hence $L$ has non-trivial cohomology on $U_y$. This proves the inclusion $\cV^{ns}(U,X)\subset \V(U,X)$.

Conversely, suppose $Rj_*(L[n])$ is simple in $Perv(X)$. That is,
$$
IC_X(L)= Rj_*(L[n]).
$$
Then, we claim that the outer terms in (\ref{eqP}) vanish simultaneously:
\begin{align*}
K(L)&:= {}^p\cH^{-1}(i_*i^{-1}Rj_*(L[n]))=0,\\
C(L)&:={}^p\cH^{0}(i_*i^{-1}Rj_*(L[n]))=0.
\end{align*}
It is clear that $C(L)$ and $K(L^{-1})$ must vanish. To show that $K(L)$ must also vanish, we recall that there is another description of $C(L)$ and $K(L)$ in terms of Deligne's nearby cycles functor. More precisely, there is an exact sequence in $Perv(X)$
$$
0\ra K(L)\ra {}^p\psi_f(L[n]){\xrightarrow{T-id}} {}^p\psi_f(L[n]) \ra C(L)\ra 0,
$$
where ${}^p\psi_f=\psi_f[-1]$ is the shifted nearby cycles functor, which restricts to a functor on perverse sheaves, and $T$ is the monodromy. 
Since the length of a perverse sheaf is additive for exact sequences, it follows that $K(L)$ and $C(L)$ have the same length as perverse sheaves. Hence $K(L)$ also vanishes.

Therefore, by (\ref{eqP}), in this case we have that $j_!(L[n])=Rj_*(L[n])$. Since $j_!L$ is the extension by zero, the stalks at points on $f^{-1}(0)$ must be zero. As before, this implies that the cohomology of $L$ in small ball complements $U_y$ of $f^{-1}(0)$ must be zero, for all $y\in f^{-1}(0)$. This shows that $\V(U,X)\subset \V^{ns}(U,X)$, and finishes the proof of the theorem.

\br\label{remBF} The proof of the inclusion (\ref{cp}) contains a small gap in \cite[5.2]{B}. It is shown there that the proof reduces to the local case. In that case, it is also shown that $\Exp(\cZ(B_F))\supset \cV^{ns}(U,X)$. This can also be seen now as an immediate consequence of Corollary \ref{corBF} of this article. The gap is in the argument that $\cV^{ns}(U,X)$ contains $\V(U,X)$. This is what we showed above, so the gap is now fixed.
\er

%%% ====================== End of main part ====================== %%%
%\newpage

%------------------------------------------------------------------
\end{document}